\documentclass[10pt, a4paper]{article} 
\usepackage{amsmath, amssymb, amsthm, amscd, ascmac, bm, enumerate, paralist, comment}
\usepackage{graphicx, xcolor, hyperref}
\usepackage[all]{xy}
\theoremstyle{definition}
	\newtheorem{Def}{Definition}
	\newtheorem{Rem}[Def]{Remark}

\theoremstyle{theorem}
	\newtheorem{Thm}[Def]{Theorem}

	\newtheorem*{MT*}{Main Theorem}

\newcommand{\R}{\mathbb{R}}
\newcommand{\C}{\mathbb{C}}
\newcommand{\Z}{\mathbb{Z}}

\renewcommand{\Re}{\operatorname{Re}}

\newcommand{\3}{{\rm{I\hspace{-.1em}I\hspace{-.1em}I}}}

\newcommand{\bash}{\backslash}

\newcommand{\ii}[2]{\langle {#1}, {#2}  \rangle}

\title{Minimal surfaces containing an epitrochoid as a geodesic}

\author{Shin Kaneda}

\date{September 16, 2022}

\begin{document}
\maketitle

\begin{abstract}
We show a complete minimal {immersion} cannot have {a certain kind of} epitrochoid as a geodesic.
\end{abstract}
\section*{Introduction}
The Bj\"{o}rling problem for minimal surfaces in $\R^3$ asks for the existence of  a minimal surface containing a given curve in $\R^3$ and a given unit normal coincides with  its Gauss map.
The representation formula of the surface of this problem's solution is given by H. A. Schwarz.
See for instance, \cite{M}, \cite{LW}, \cite{MW}, and so on.

\begin{figure}[htbp]
\centering
	\begin{tabular}{c}
		\begin{minipage}{0.33\hsize}
		\centering
		\includegraphics[scale=0.15]{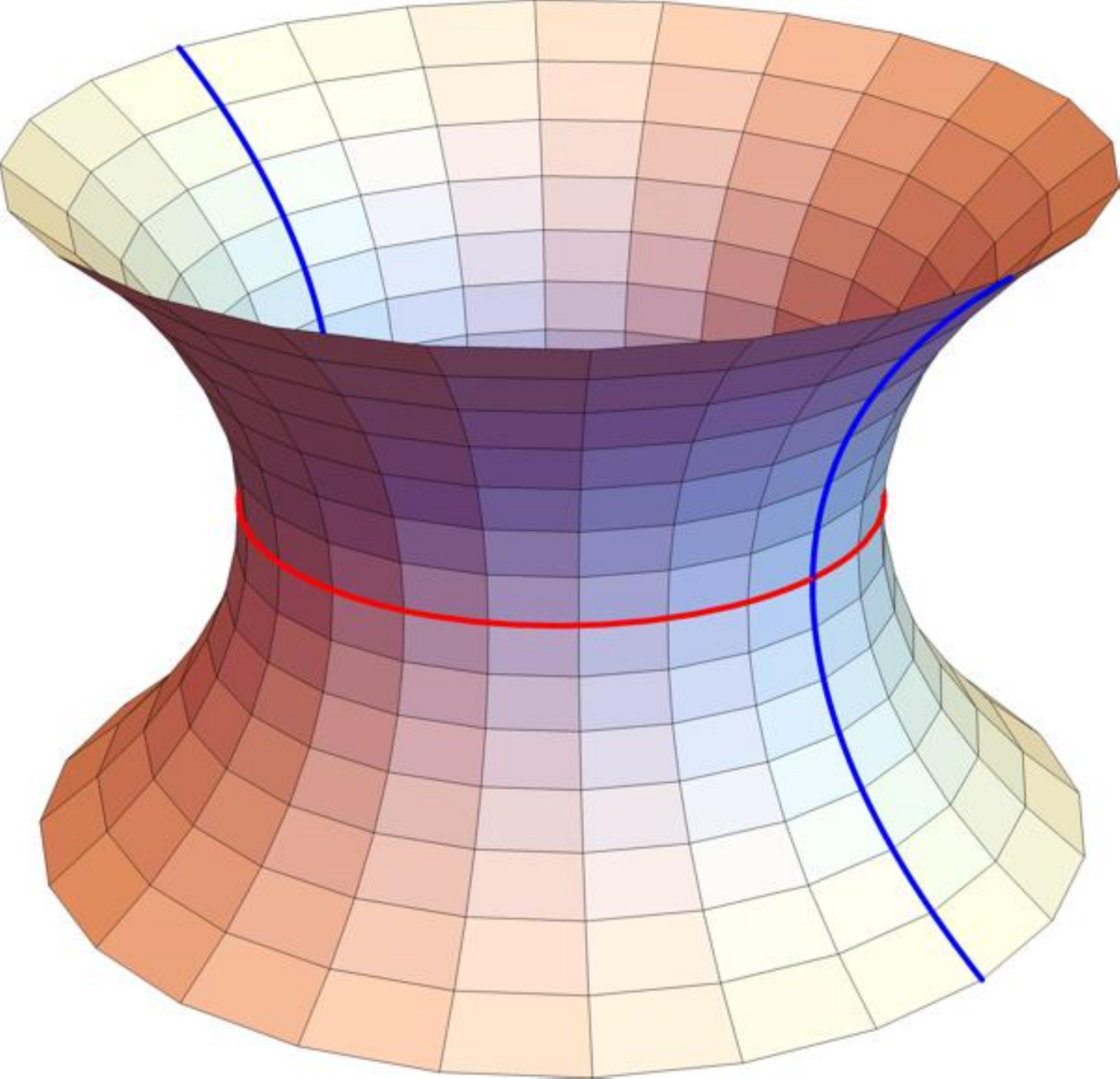}
		\end{minipage}
		
		\begin{minipage}{0.33\hsize}
		\centering
		\includegraphics[scale=0.15]{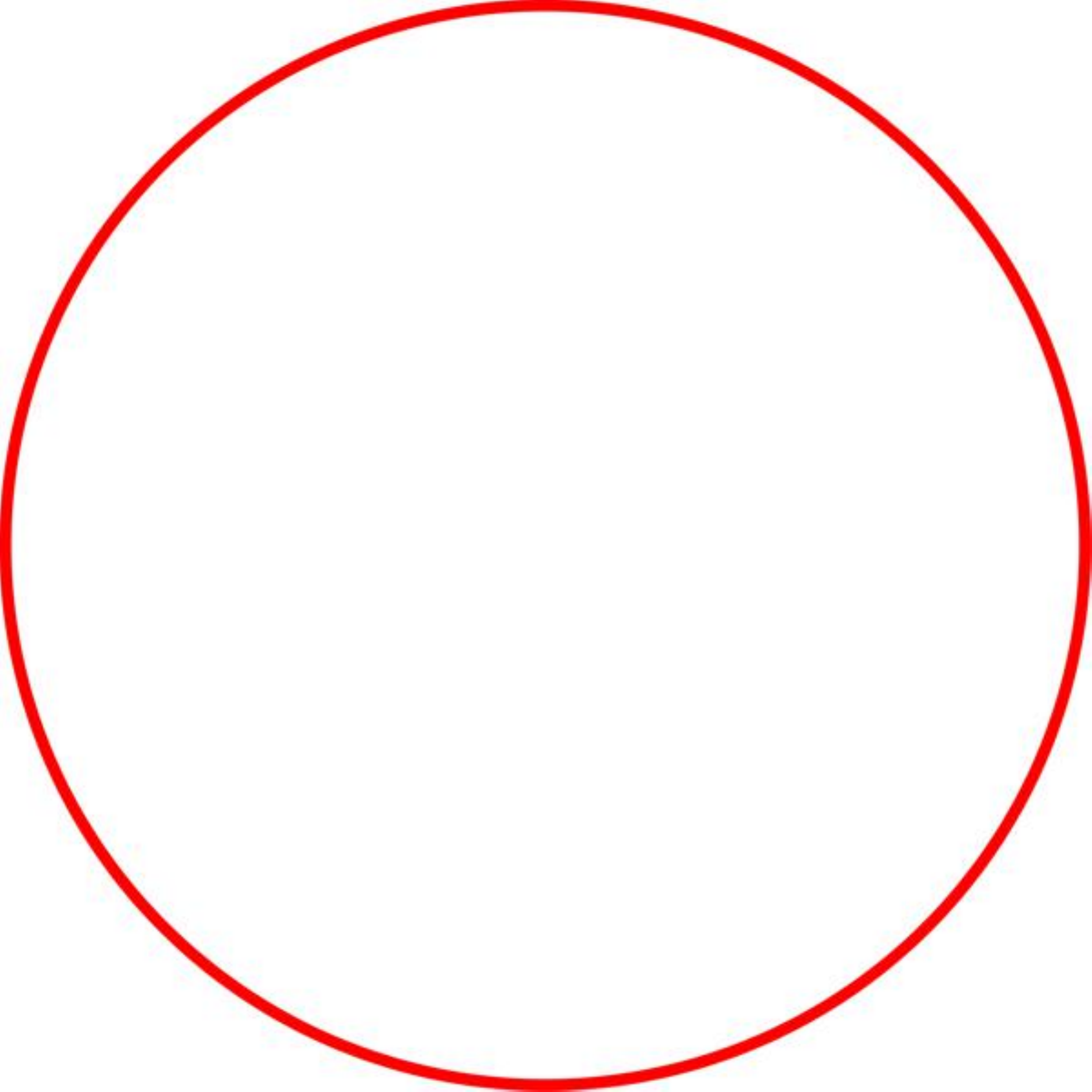}
		\end{minipage}
				
		\begin{minipage}{0.33\hsize}
		\centering
		\includegraphics[scale=0.15]{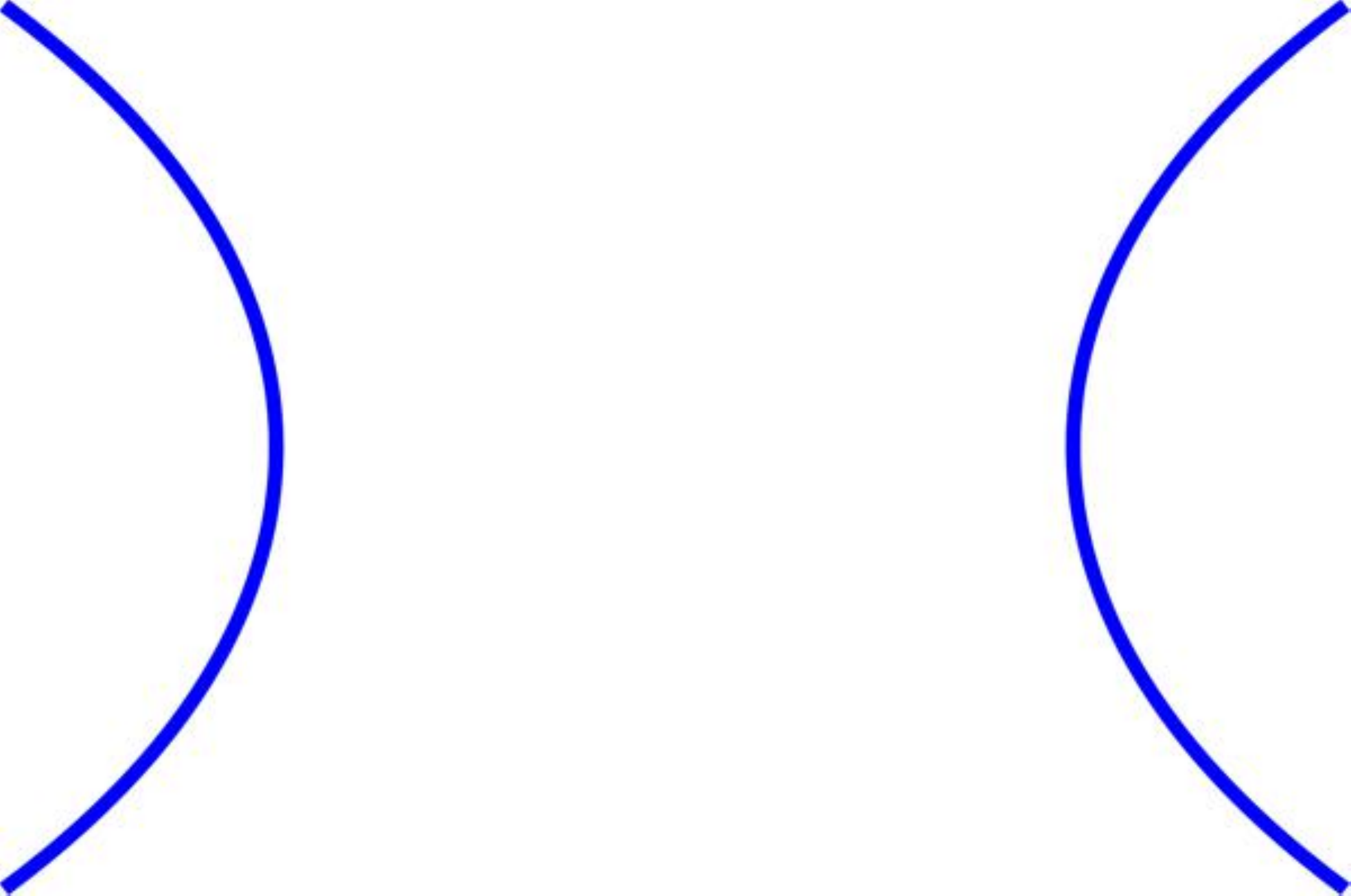}
		\end{minipage}
	\end{tabular}
\caption{
Catenoid and planar curves, a circle and catenaries, contained as geodesics.
}\label{fig:cat}
\end{figure}

\begin{figure}[htbp]
\centering
	\begin{tabular}{c}
		\begin{minipage}{0.33\hsize}
		\centering
		\includegraphics[scale=0.15]{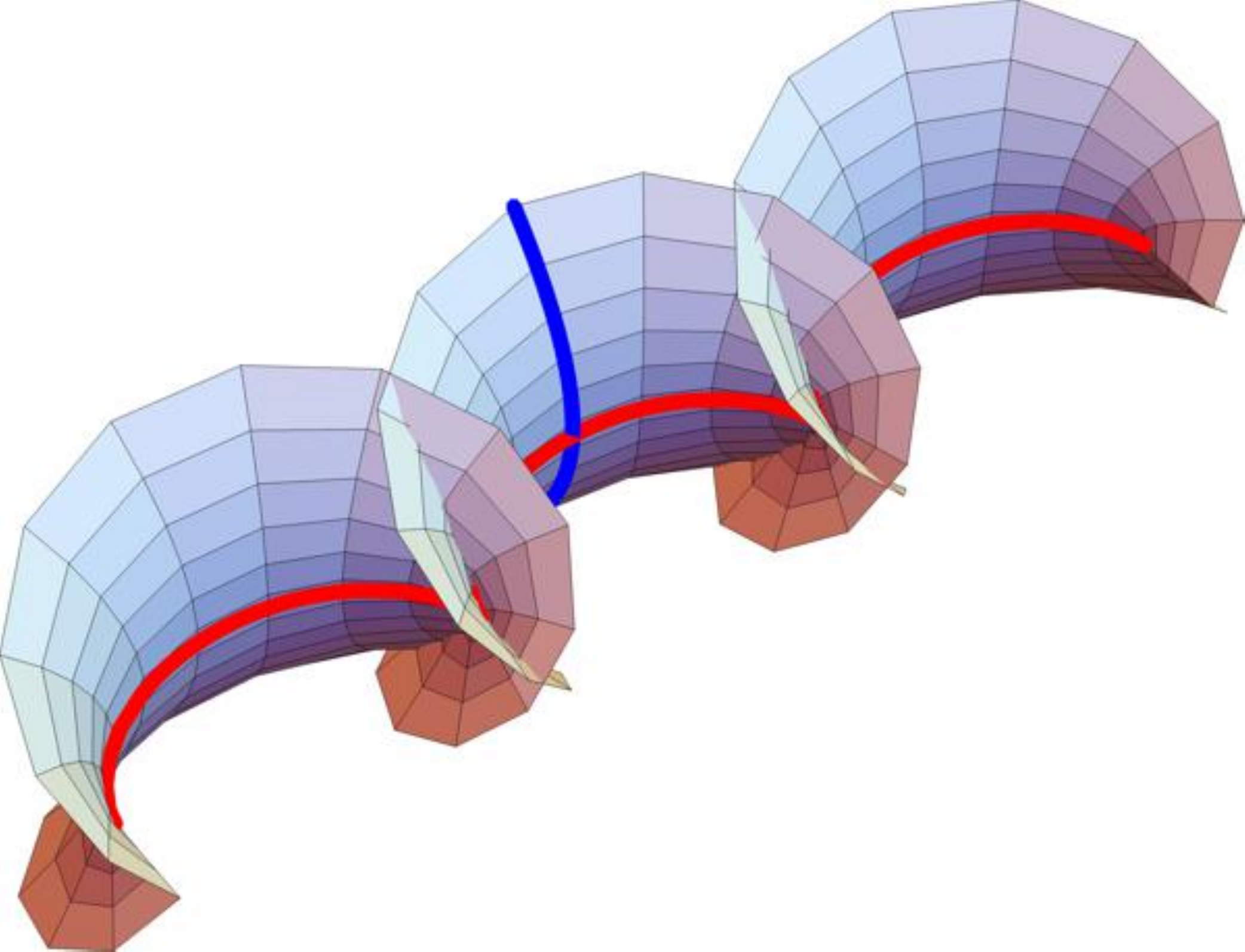}
		\end{minipage}
		
		\begin{minipage}{0.33\hsize}
		\centering
		\includegraphics[scale=0.15]{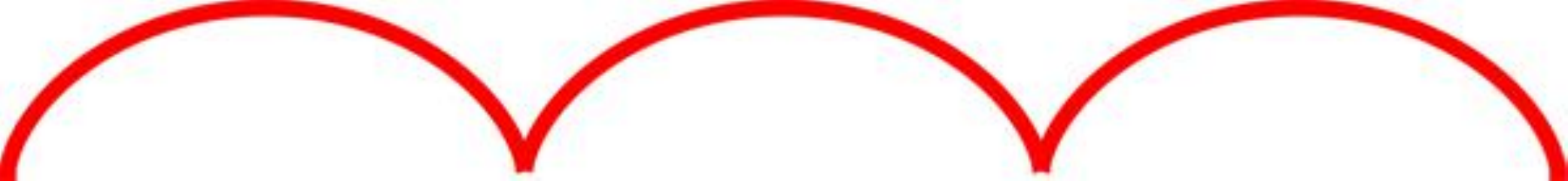}
		\end{minipage}
				
		\begin{minipage}{0.33\hsize}
		\centering
		\includegraphics[scale=0.15]{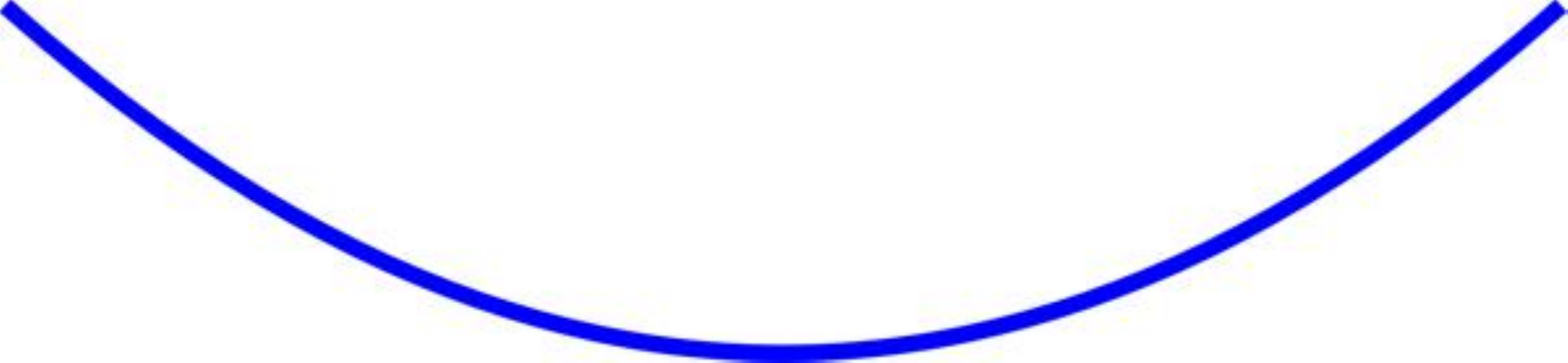}
		\end{minipage}
	\end{tabular}
\caption{
Catalan surface and planar curves, a cycloid and a parabola, contained as geodesics.
}\label{fig:catalan}
\end{figure}

We consider minimal surfaces containing planar curves as geodesics.
For example, a catenoid contains a circle and catenaries as geodesics (See Figure \ref{fig:cat}), Catalan surface contains a cycloid and a parabola (See Figure \ref{fig:catalan}). 
The minimal surfaces containing a conic curve as a geodesic are studied in \cite{AGS}, \cite{P}, and so on.
Note that the Catalan surface and the minimal surfaces containing a conic curve as a geodesic are not {immersed}, while the catenoid is {immersed and} complete.

\begin{figure}[htbp]
\centering
	\begin{tabular}{c}
		\begin{minipage}{0.33\hsize}
		\centering
		\includegraphics[scale=0.15]{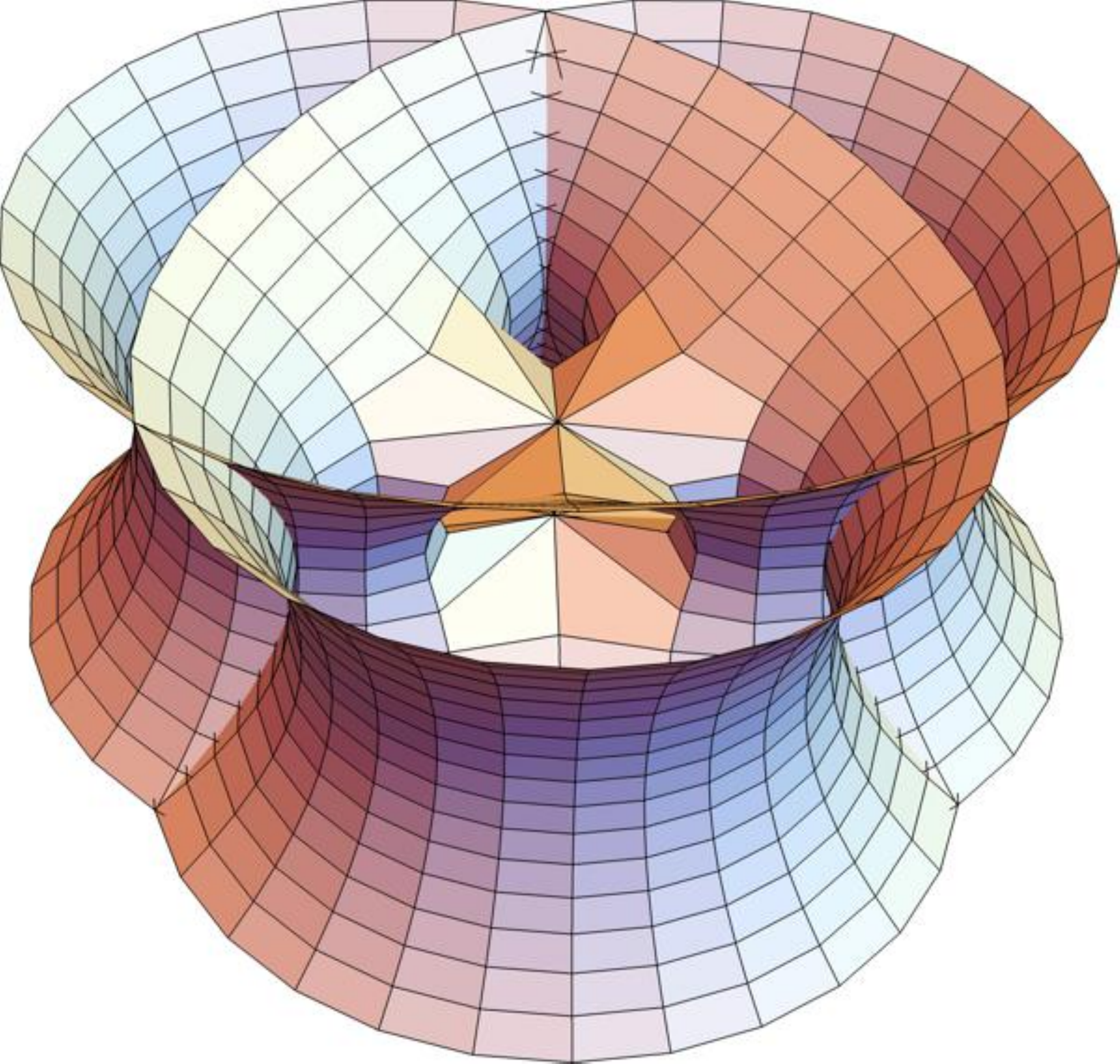}
		\end{minipage}
		
		\begin{minipage}{0.33\hsize}
		\centering
		\includegraphics[scale=0.15]{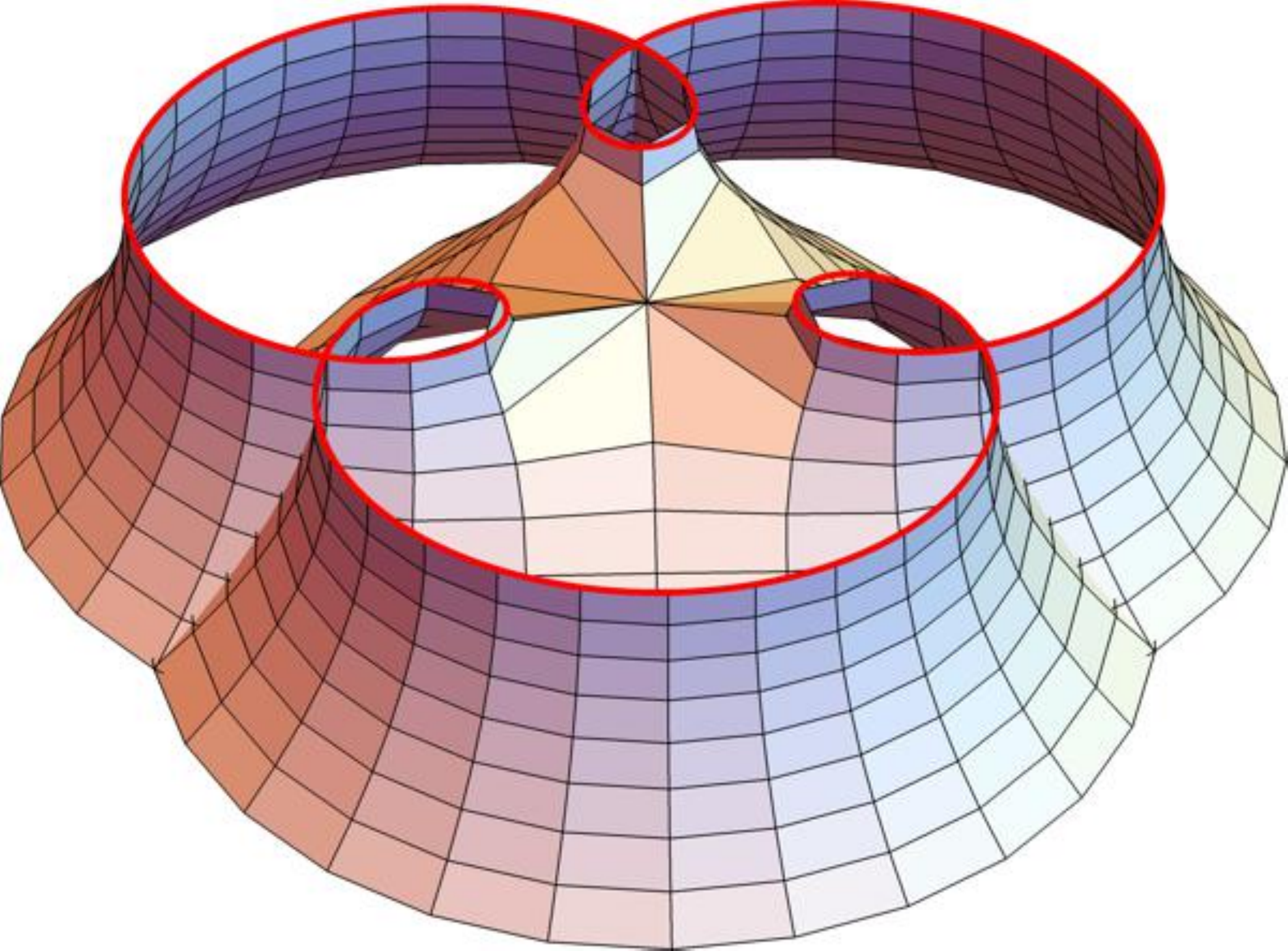}
		\end{minipage}
				
		\begin{minipage}{0.33\hsize}
		\centering
		\includegraphics[scale=0.15]{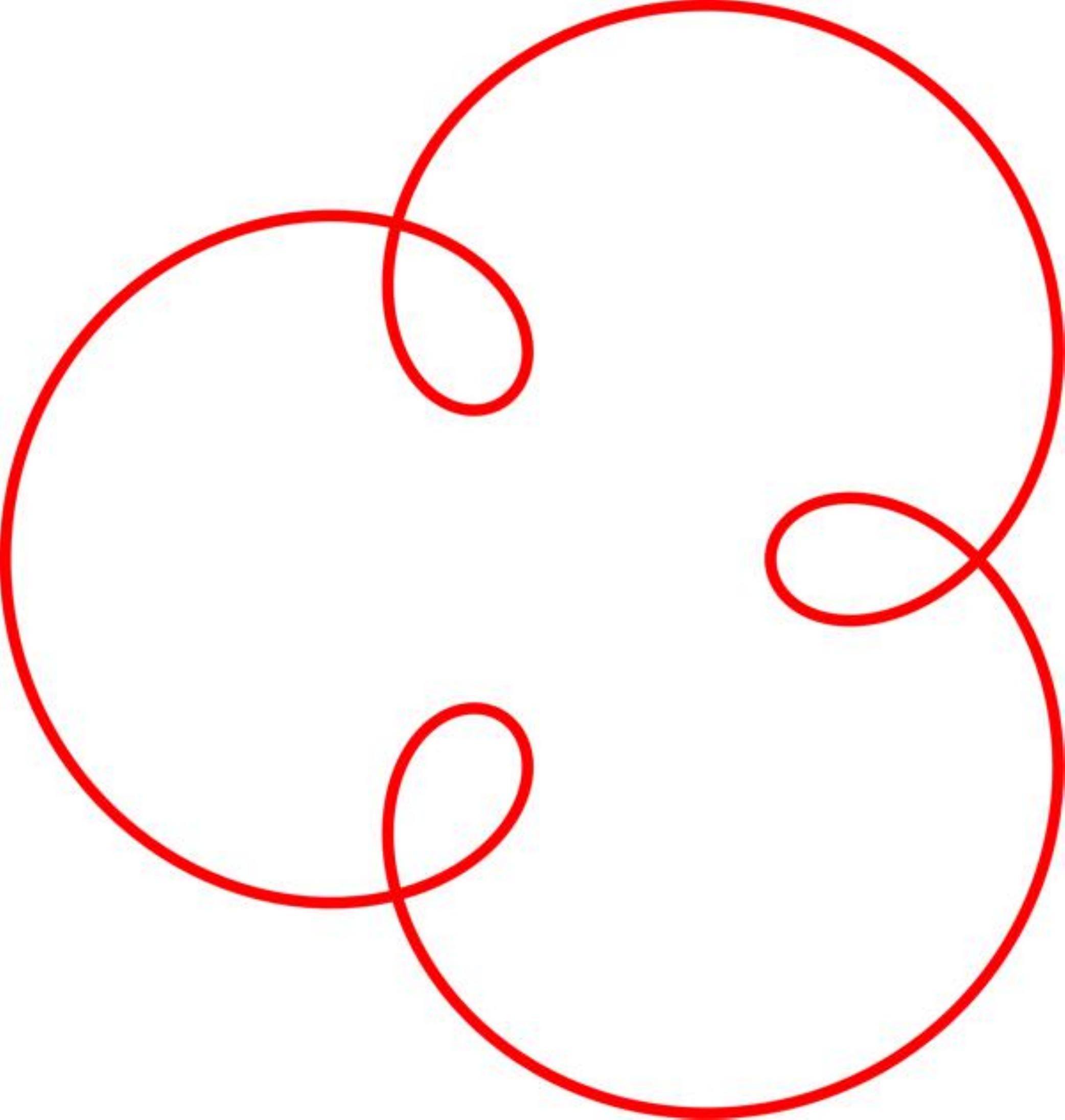}
		\end{minipage}
	\end{tabular}
\caption{The minimal surface constructed in \cite{FS} (left), a half cut away from the surface by $x_1x_2$-plane (middle), and  its intersection with the $x_1x_2$-plane (right).
}\label{fig:fs}
\end{figure}

In \cite{FS}, Fujimori and Shoda constructed complete minimal {immersions} with finite total curvature.
These surfaces contain a planar curve as a geodesic.
These curves really look like  epitrochoids (See Figure \ref{fig:fs}).
So the author tried to prove that these curves are indeed epitrochoids.
Then, surprisingly, we found that these curves are different from the epitrochoids.
In fact, we prove the following : 
\begin{MT*}
A complete minimal {immersion} in $\R^3$ cannot have {a single-wrapped} epitrochoid as a geodesic.
\end{MT*}
{For the definition of single-wrapped epitrochoid, see Definition \ref{def:single}.}
By this theorem, we see that the surfaces constructed in \cite{FS} do not contain an epitrochoid as a geodesic because these surfaces are {immersed and} complete.

\section{Preliminaries}
\begin{Thm}[The Weierstrass representation {\cite[Theorem 22.27]{AGS}}]\label{th:wr}
Let $(g, \eta)$ be a pair of a meromorphic function $g$ and a holomorphic differential $\eta$ on a Riemann surface $M$ such that 
\begin{equation}\label{eq:complete}
(1+|g|^2)^2\eta\bar{\eta}
\end{equation}
gives a Riemannian metric on $M$. We set
\begin{equation}\label{eq:wr}
\Phi =
	\begin{pmatrix}
	\phi_1 \\
	\phi_2 \\
	\phi_3
	\end{pmatrix}
=
	\begin{pmatrix}
	(1 - g^2)\eta \\
	i(1 + g^2)\eta \\
	2g\eta
	\end{pmatrix}
. 
\end{equation}
Suppose that
\begin{equation}\label{eq:period}
\Re \oint_{\gamma} \Phi = {\bf{0}}
\end{equation}
holds for any $\gamma \in H_1(M, \Z).$ Then
\begin{equation}
f(z) = \Re \int_{z_0}^{z} \Phi : M \to \R^3, \qquad (z_0 \in M)
\end{equation}
defines a minimal {immersion} {in $\R^3$}. 
\end{Thm}
The pair $(g, \eta)$ in Theorem \ref{th:wr} is called the Weierstrass data of $f$.
\begin{Rem}
Let $(g, \eta)$ be the Weierstrass data of $f$. Then $(g, \eta)$ can be written as
\begin{equation}\label{eq:geta}
g = \frac{\phi_3}{\phi_1-i\phi_2}, \qquad \eta = \phi_1 - i \phi_2.
\end{equation}
Moreover, $g$ coincides with the composition of the unit normal $\nu$ of $f$ and the stereographic projection $\sigma$ from the north pole, that is, $g = \sigma \circ \nu$.
Thus we call $g$ the Gauss map of $f$.
\end{Rem}
	
\begin{Thm}[Schwarz {\cite[Section 22.6]{AGS}}]
Let $c : I \to \R^3$ be a real analytic regular curve and $n : I \to \R^3$  a real analytic vector field along $c$ such that
\[
\ii{n(t)}{n(t)}=1, \qquad \ii{n(t)}{c'(t)}=0, \qquad t \in I,
\]
where $I$ is an interval in $\R$.
Then, for sufficiently small $\varepsilon > 0$ there exists a unique minimal immersion $f : I \times (-\varepsilon, \varepsilon) \to \R^3$ satisfying
\[
f(t) = c(t), \qquad \nu(t) = n(t), \qquad t \in I,
\]
where $\nu$ is the unit normal of $f$.
Moreover, $f$ can be written on $I \times (-\varepsilon, \varepsilon)$ as
\begin{equation}\label{eq:sch}
f(z) = \Re \left\{ c(z) -i \int_{z_0}^z n(w) \times c'(w) dw\right\},
\end{equation}
where $c(z), n(z)$ are extensions of $c(t), n(t)$ to functions of a complex variable $z = t + i s \in I \times (-\varepsilon, \varepsilon)$.
\end{Thm}
	
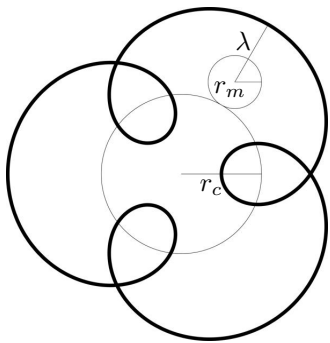
\begin{figure}[htbp]
	\centering
	\includegraphics[scale=0.2]{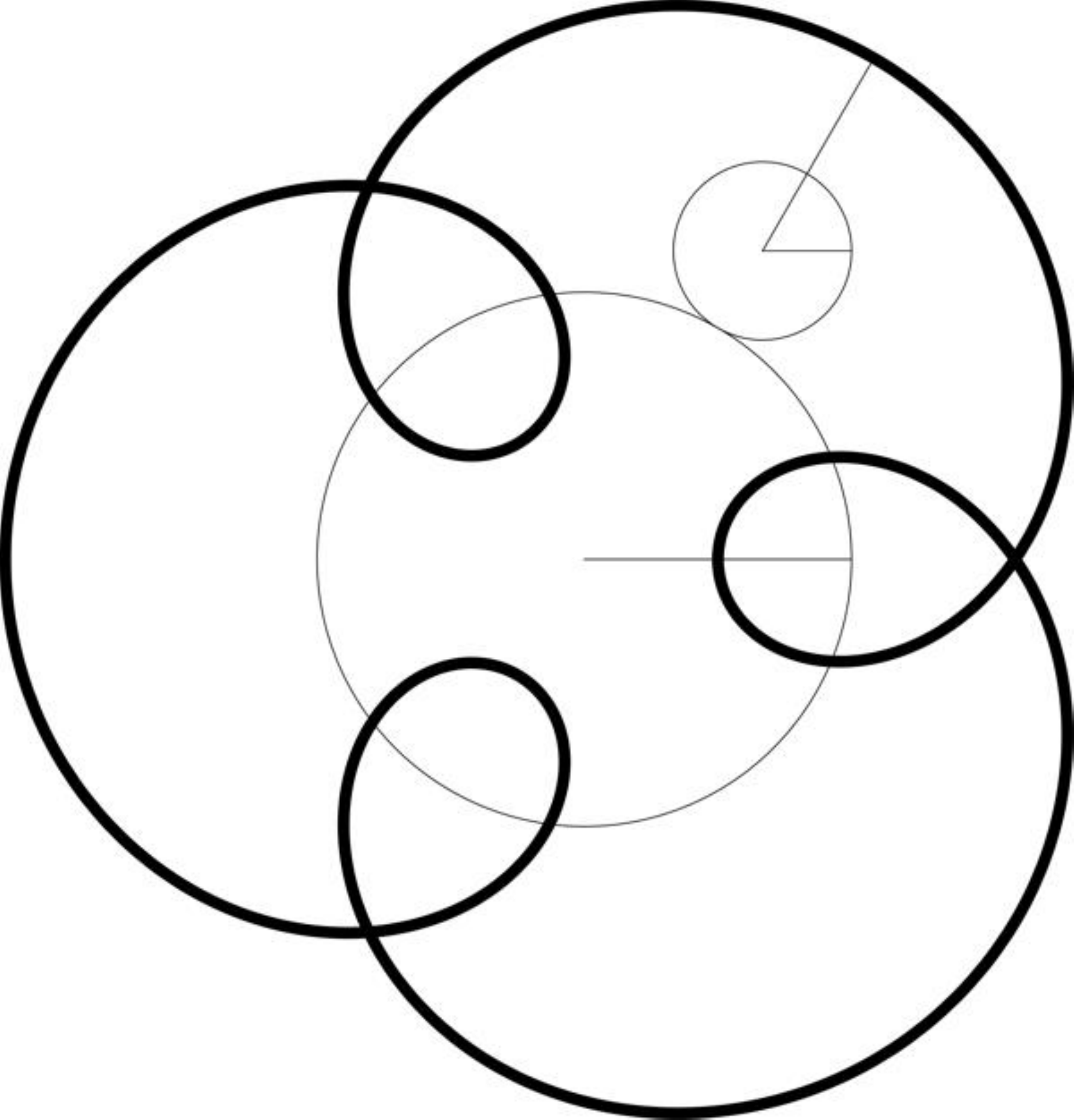}
	\unitlength=1pt
	\begin{picture}(0,0)
		\put(-47,58){\makebox(0,0)[cc]{$r_c$}}
		\put(-40,95){\makebox(0,0)[cc]{$r_{m}$}}
		\put(-35,113){\makebox(0,0)[cc]{$\lambda$}}
	\end{picture}
	\caption{The thick curve indicates an epitrochoid, the big  circle is a fixed circle, and the small  circle is a rolling circle.}\label{fig:howtodrawepi}
\end{figure}

We consider minimal {immersions} containing a real analytic regular planar curve that is not a straight line $c(t)=(x(t), y(t), 0)$ as a geodesic.
We also suppose that $c(t)$ is regular, that is, $x'(t)^2+y'(t)^2>0$ for any $t \in I$.
Since $c''(t)$ is contained in the $x_1x_2$-plane  and $c$ is a geodesic of $f$, the unit normal $n$ of $c$ can be written as
\[ 
n(t) = \frac{1}{\sqrt{x'(t)^2+y'(t)^2}}(-y'(t), x(t), 0).
\] 
Applying (\ref{eq:sch}) to $c$ and $n$, we have
\begin{equation}\label{eq:scpl}
f(z) = \Re \int_{z_0}^{z}\left(x'(w), y'(w), i\sqrt{x'(w)^2+y'(w)^2}\right)dw.
\end{equation}
Moreover, from (\ref{eq:geta}), we have
\begin{equation}\label{eq:bjwr}
g = i \sqrt{\frac{x'(z)+iy'(z)}{x'(z)-iy'(z)}}, \qquad \eta = (x'(z)-iy'(z))dz.
\end{equation}
Let $c(t)=(x(t), y(t),0)$ be an epitrochoid.  
Then it is written by
\begin{equation}\label{eq:epi1}
	\begin{cases}
	x(t)=(r_c + r_m)\cos t-\lambda\cos\left(\dfrac{r_c+r_m}{r_m}t\right),\vspace{4pt}\\
	y(t)=(r_c + r_m)\sin t-\lambda\sin\left(\dfrac{r_c+r_m}{r_m}t\right),
	\end{cases}
\end{equation}
where $r_c>0$ and $r_m>0$ are radii of a fixed circle and a rolling circle around outside of a fixed circle respectively, $\lambda>0$ is a distance between a center of a rolling circle and a point on the epitrochoid (See Figure \ref{fig:howtodrawepi}).
For the curve $c(t)$ to be the regular, we assume that $\lambda \neq r_m$.

{Without loss of generality, we may assume that $r_c = 1$}.
\begin{Def}\label{def:single}
{An epitrochoid (\ref{eq:epi1}) is called {\it{single-wrapped}} if there exists $k\in\Z_{>0}$ such that $r_m=1/(k+1)$.
Note that if an epitrochoid is single-wrapped, then it is closed when the rolling circle wraps once around the fixed circle.}
\end{Def}
{We consider a single-wrapped epitrochoid, that is, we set} $r_m = 1/(k+1), k\in\Z_{>0}$. 
{Then} (\ref{eq:epi1}) becomes
\[
	\begin{cases}
	x(t)=\left(1 + \dfrac{1}{k+1}\right)\cos t-\lambda\cos\left(\dfrac{1+1/(k+1)}{1/(k+1)}t\right),\vspace{4pt}\\
	y(t)=\left(1 + \dfrac{1}{k+1}\right)\sin t-\lambda\sin\left(\dfrac{1+1/(k+1)}{1/(k+1)}t\right).
	\end{cases}
\]

\begin{figure}[htbp]
\centering
	\begin{tabular}{c}
		\begin{minipage}{0.33\hsize}
		\centering
		\includegraphics[scale=0.15]{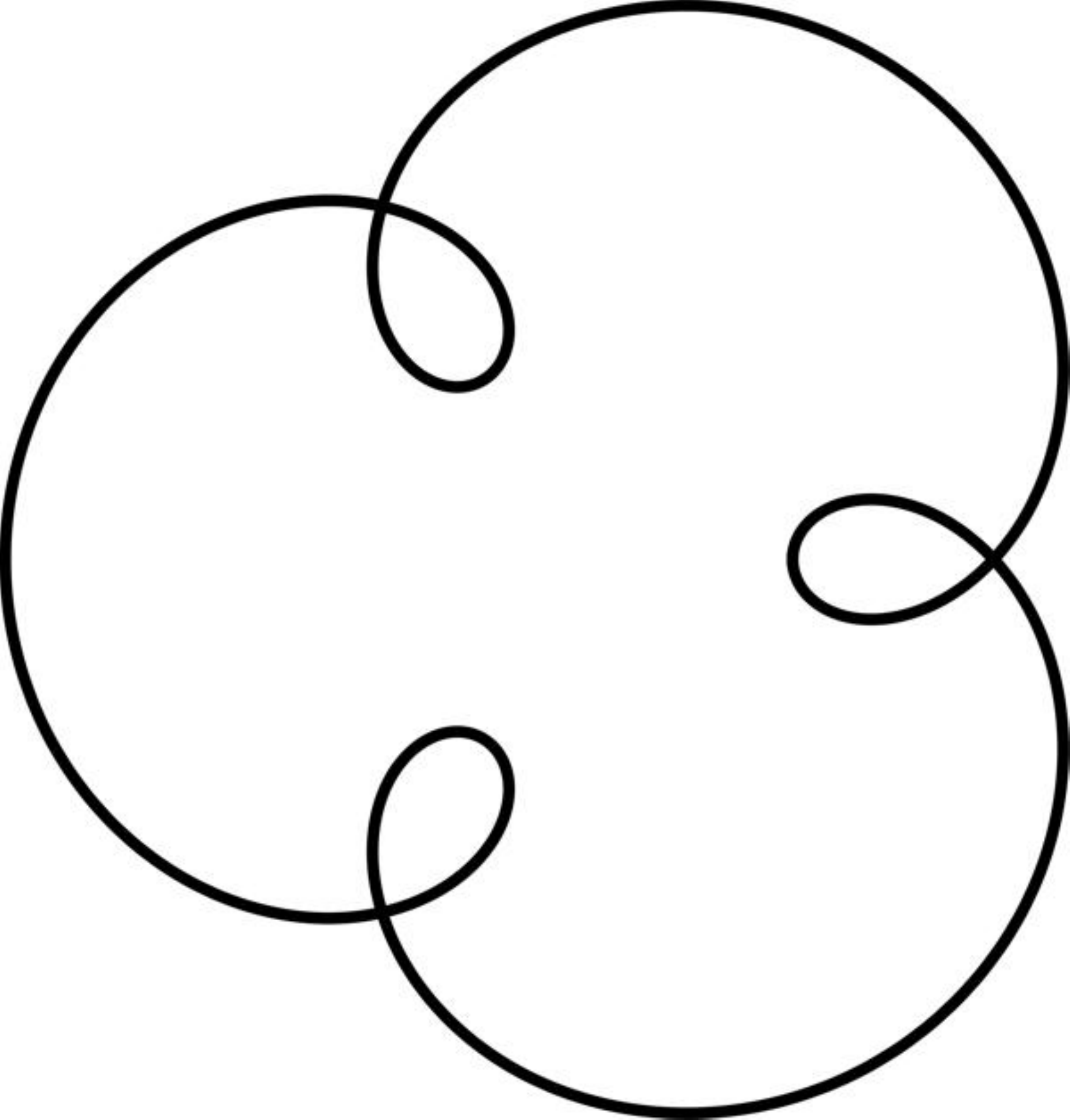}
		\end{minipage}
		
		\begin{minipage}{0.33\hsize}
		\centering
		\includegraphics[scale=0.15]{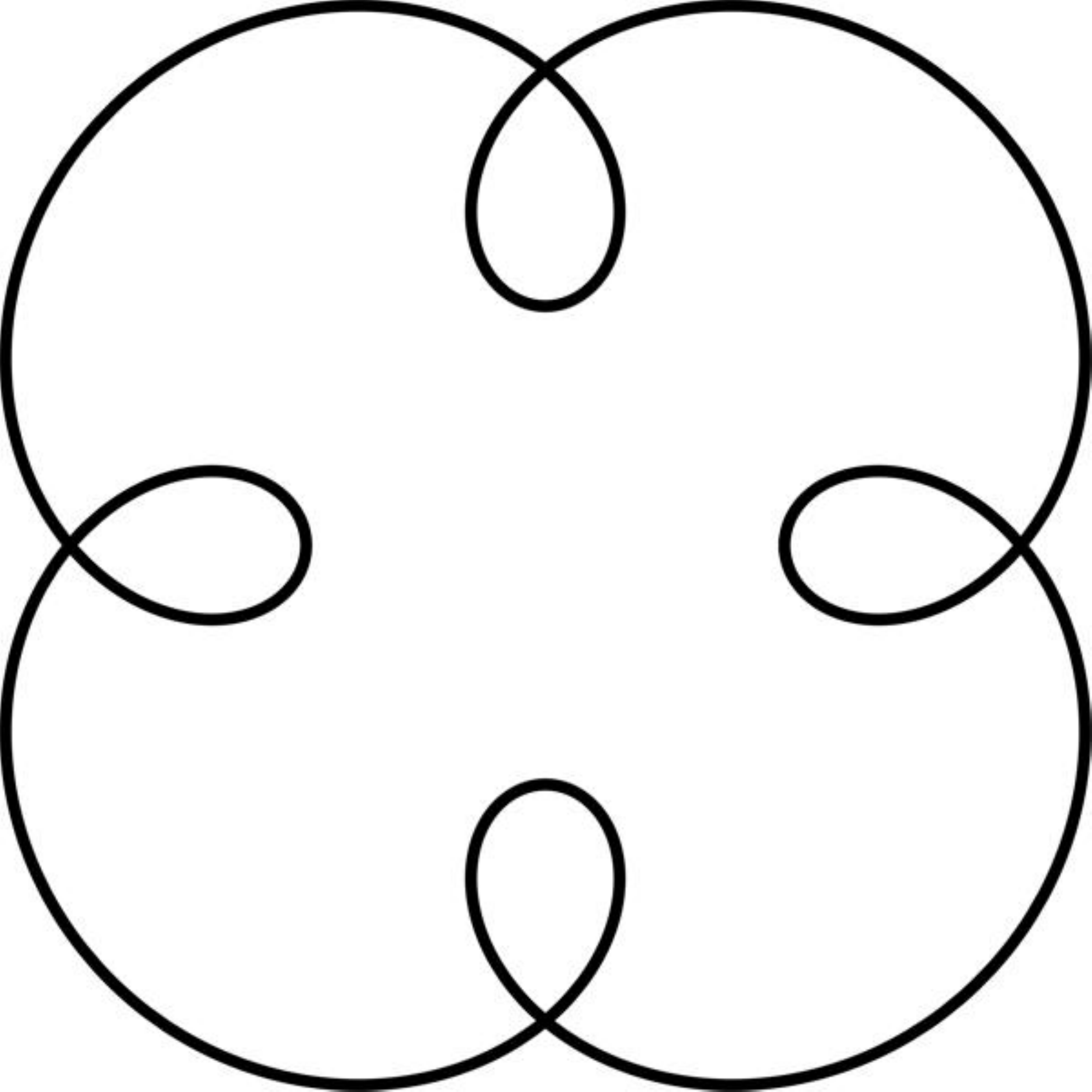}
		\end{minipage}
				
		\begin{minipage}{0.33\hsize}
		\centering
		\includegraphics[scale=0.15]{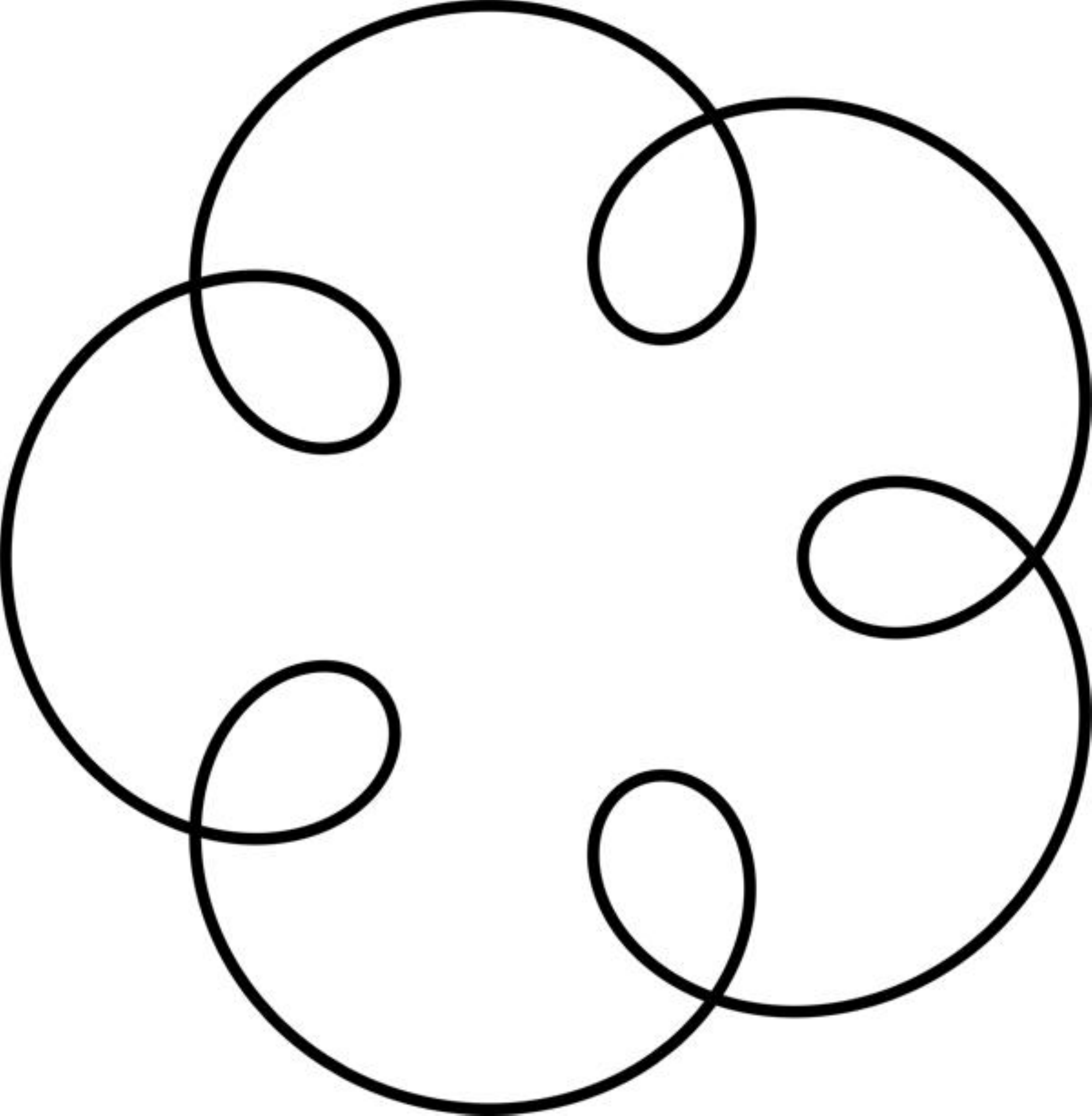}
		\end{minipage}
	\end{tabular}
\caption{{Single-wrapped epitrochoids} (left : $k=2$, middle : $k=3$, right : $k=4$) (${\lambda} > r_m$)}\label{fig:epi}
\end{figure}

Scaling up this curve, we can set
\begin{equation}\label{eq:epi2}
	\begin{cases}
	x(t)=\left(k+2\right)\cos t-(k+1)\lambda\cos\left(k+2\right)t\vspace{4pt}\\
	y(t)=\left(k+2\right)\sin t-(k+1)\lambda\sin\left(k+2\right)t.
	\end{cases}
\end{equation}
See Figure \ref{fig:epi}.

\section{The proof of the main theorem}
We consider a {complete} minimal {immersion} containing {a single-wrapped} epitrochoid as a geodesic.
Applying (\ref{eq:scpl}) and (\ref{eq:epi2}) , we have
\begin{equation}
f(z) = \Re \int_{z_0}^{z}
	\begin{pmatrix}
	-(k+2)\sin z + (k+1)(k+2)\lambda \sin(k+2)z \\
	(k+2)\cos z -(k+1)(k+2)\lambda \cos(k+2) z \\
	i(k+2) \sqrt{1+(k+1)^2\lambda^2-2\lambda(k+1)\cos(k+1)z}
	\end{pmatrix}
dz.
\end{equation}
By (\ref{eq:bjwr}), 
\[
g = i \sqrt{\frac{e^{iz}-\lambda(k+1)e^{(k+2)iz}}{-e^{-iz}+\lambda(k+1)e^{-(k+2)iz}}}, \qquad \eta = (k+2)(e^{-iz}+\lambda(k+1)e^{-(k+2)iz})dz.
\]
After the substitution $v = e^{iz}$ and rotating the surface in $\R^3$ by $-\pi/2$ about the $x_3$-axis, we have
\[
g = \sqrt{\frac{v^{k+3}(1-\lambda(k+1)v^{k+1})}{-v^{k+1}+\lambda(k+1)}}, \qquad \eta = (k+2)\frac{v^{k+1}-\lambda(k+1)}{v^{k+3}}dv.
\]
We set 
\[
w^2 = 
	\begin{cases}
	v(1-\lambda(k+1)v^{k+1})(\lambda(k+1)-v^{k+1}) \qquad \text{ ($k$ is even),}\\
	(1-\lambda(k+1)v^{k+1})(\lambda(k+1)-v^{k+1}) \qquad \text{ ($k$ is odd).}
	\end{cases}
\]

\noindent{\bf{Case 1.}} $k$ is even. \vspace{3.0mm} 

Let $\overline{M}$ be the compact Riemann surface of genus $k+1$ defined by 
\[
\overline{M} = \{(v, w) \in (\C\cup\{\infty\})^2 \: |\: w^2 = v(1-\lambda(k+1)v^{k+1})(\lambda(k+1)-v^{k+1})\},
\]
and set $M = \overline{M}\bash\{(0,0), (\infty,\infty)\}.$
Then the Weierstrass data of $f$ can be written as
\[
g = -\frac{w v^{\frac{k+2}{2}}}{v^{k+1}-\lambda(k+1)}, \qquad \eta = \frac{v^{k+1}-\lambda(k+1)}{v^{k+1}}dv.
\]
{It is easy to check that  (\ref{eq:complete}) degenerates at $v = \sqrt[k+1]{\lambda(k+1)}$ and $v = \sqrt[k+1]{1/(\lambda(k+1))}$.
Therefore $f$ is not an immersion at these points.} 
See Table \ref{tab:data-even}.
\begin{table}[h]
\centering
	\begin{tabular}{c|cccccc}
	$(v,w)$&$(0,0)$&$(\sqrt[k+1]{\lambda(k+1)}, 0)$&$(\sqrt[k+1]{1/(\lambda(k+1))}, 0)$&$(\infty, \infty)$ \\ \hline 
	$g$&$0^{k+3}$&$\infty^1$&$0^1$&$\infty^{k+3}$\\ 
	$\eta$&$\infty^{2k+5}$&$0^{3}$&$0^1$&$0^1$
	\end{tabular}
\caption{Orders of zeros and poles of $g$ and $\eta$}\label{tab:data-even}
\end{table}

\begin{figure}[htbp]
\centering
	\begin{tabular}{c}
		\begin{minipage}{0.3\hsize}
		\centering
		\includegraphics[scale=0.15]{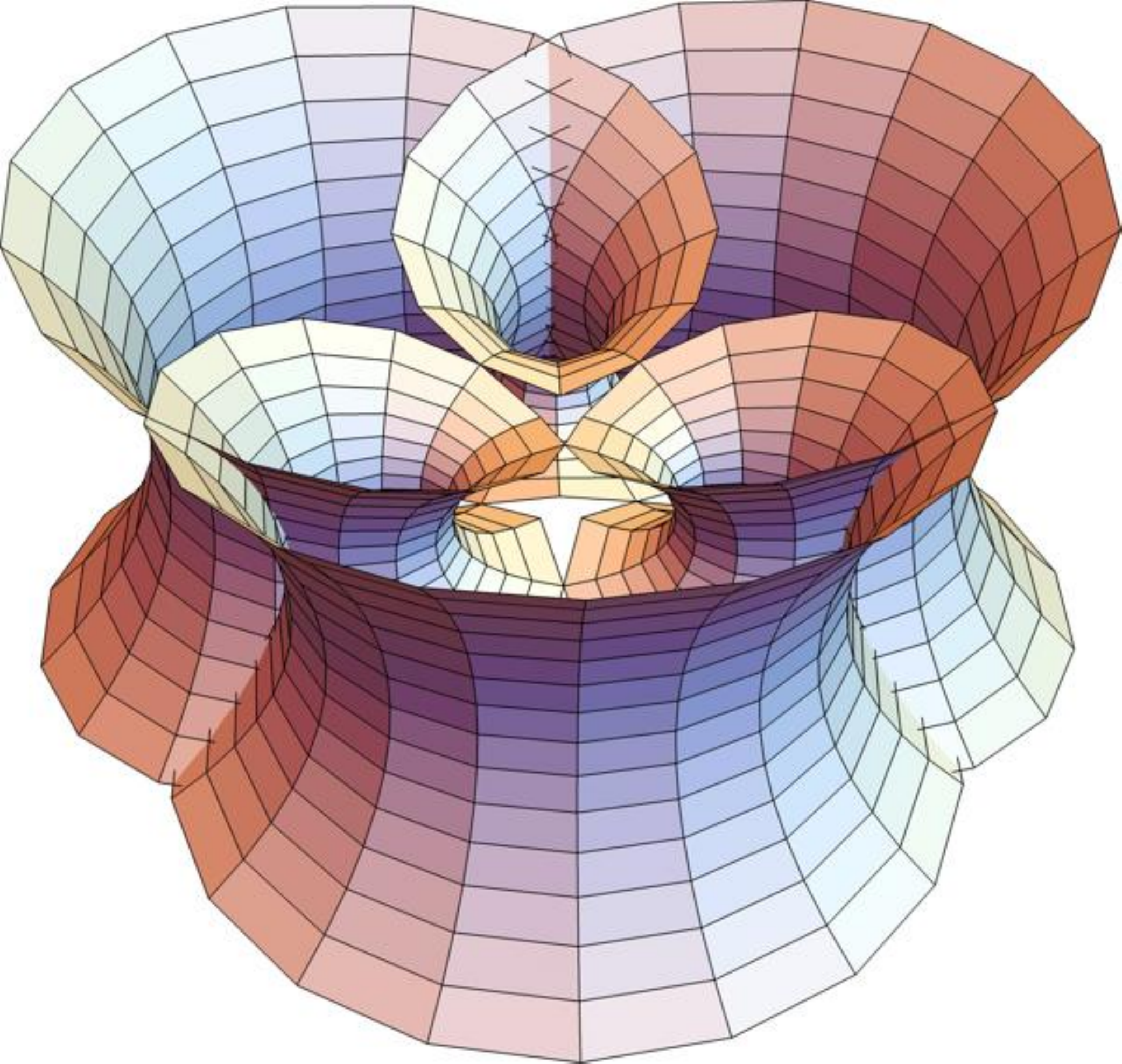}
		\end{minipage}
		
		\begin{minipage}{0.3\hsize}
		\centering
		\includegraphics[scale=0.15]{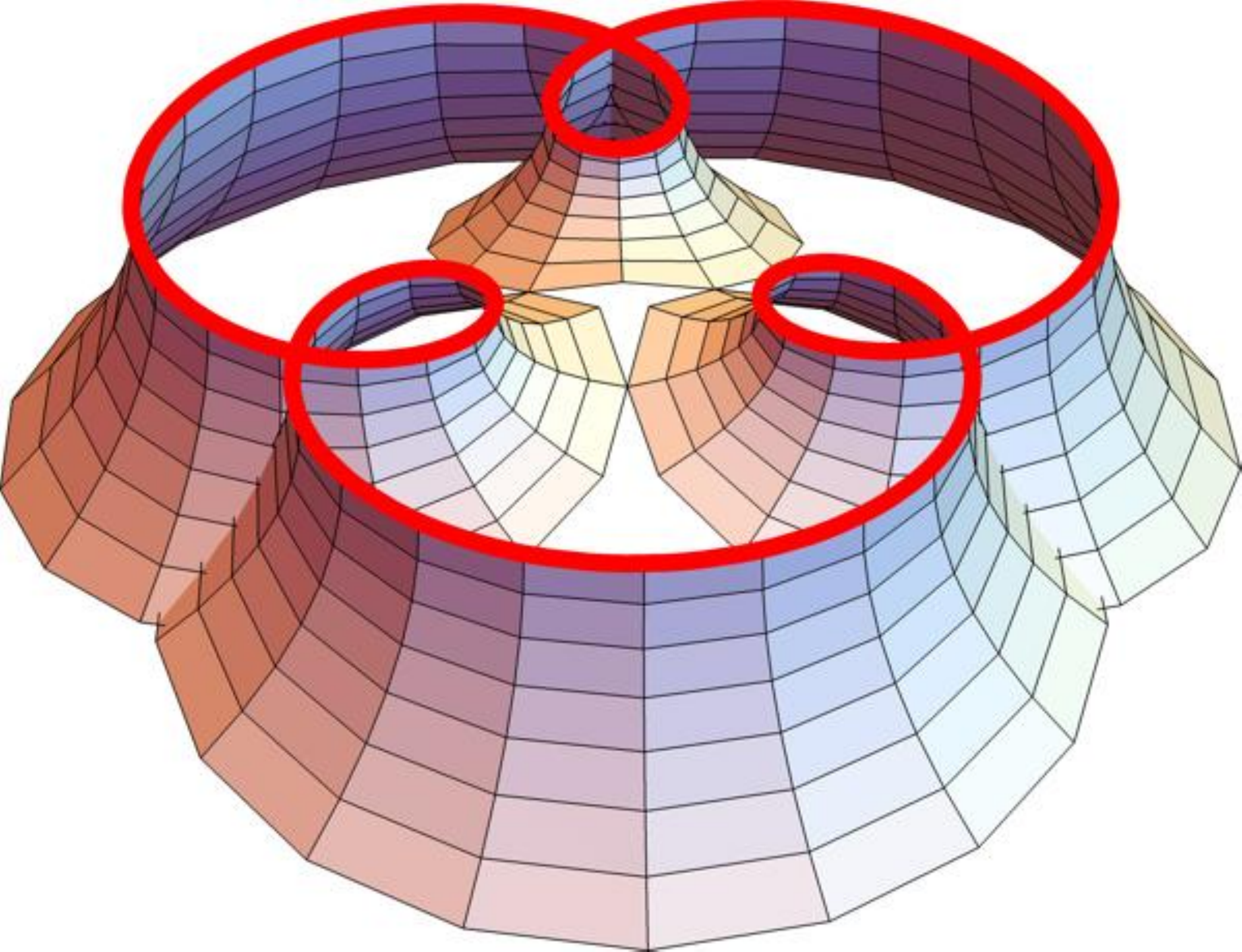}
		\end{minipage}
	\end{tabular}
\caption{
The surface with $k=2$ (left), and a half cut away from the surface by the $x_1x_2$-plane (right).
}
\end{figure}
$\\$
\noindent{\bf{Case 2.}} $k$ is odd. \vspace{3.0mm} 

Let $\overline{M}$ be the compact Riemann surface of genus $k$ defined by 
\[
\overline{M} = \{(v, w) \in (\C\cup\{\infty\})^2 \: |\: w^2 = (1-\lambda(k+1)v^{k+1})(\lambda(k+1)-v^{k+1})\},
\]
and set $M = \overline{M}\bash\{(0,\pm{\sqrt{\lambda(k+1)}}), (\infty,\infty)\}.$
Then the Weierstrass data of $f$ can be written as
\[
g = -\frac{w v^{\frac{k+3}{2}}}{v^{k+1}-\lambda(k+1)}, \qquad \eta = \frac{v^{k+1}-\lambda(k+1)}{v^{k+1}}dz.
\]
{It is easy to check that  (\ref{eq:complete}) degenerates at $v = \sqrt[k+1]{\lambda(k+1)}$ and $v = \sqrt[k+1]{1/(\lambda(k+1))}$.
Therefore $f$ is not an immersion at these points.} 
See Table \ref{tab:data-odd}.
\begin{table}[h]
\centering
	\begin{tabular}{c|cccccc}
	$(v,w)$&$(0,\pm{\sqrt{\lambda(k+1)}})$&$(\sqrt[k+1]{\lambda(k+1)}, 0)$&$(\sqrt[k+1]{1/(\lambda(k+1))}, 0)$&$(\infty, \infty)$ \\ \hline 
	$g$&$0^{\frac{k+3}{2}}$&$\infty^1$&$0^1$&$\infty^{\frac{k+3}{2}}$\\ 
	$\eta$&$\infty^{k+3}$&$0^{3}$&$0^1$&$-$
	\end{tabular}
\caption{Orders of zeros and poles of $g$ and $\eta$}\label{tab:data-odd}
\end{table}

\begin{figure}[htbp]
\centering
	\begin{tabular}{c}
		\begin{minipage}{0.3\hsize}
		\centering
		\includegraphics[scale=0.15]{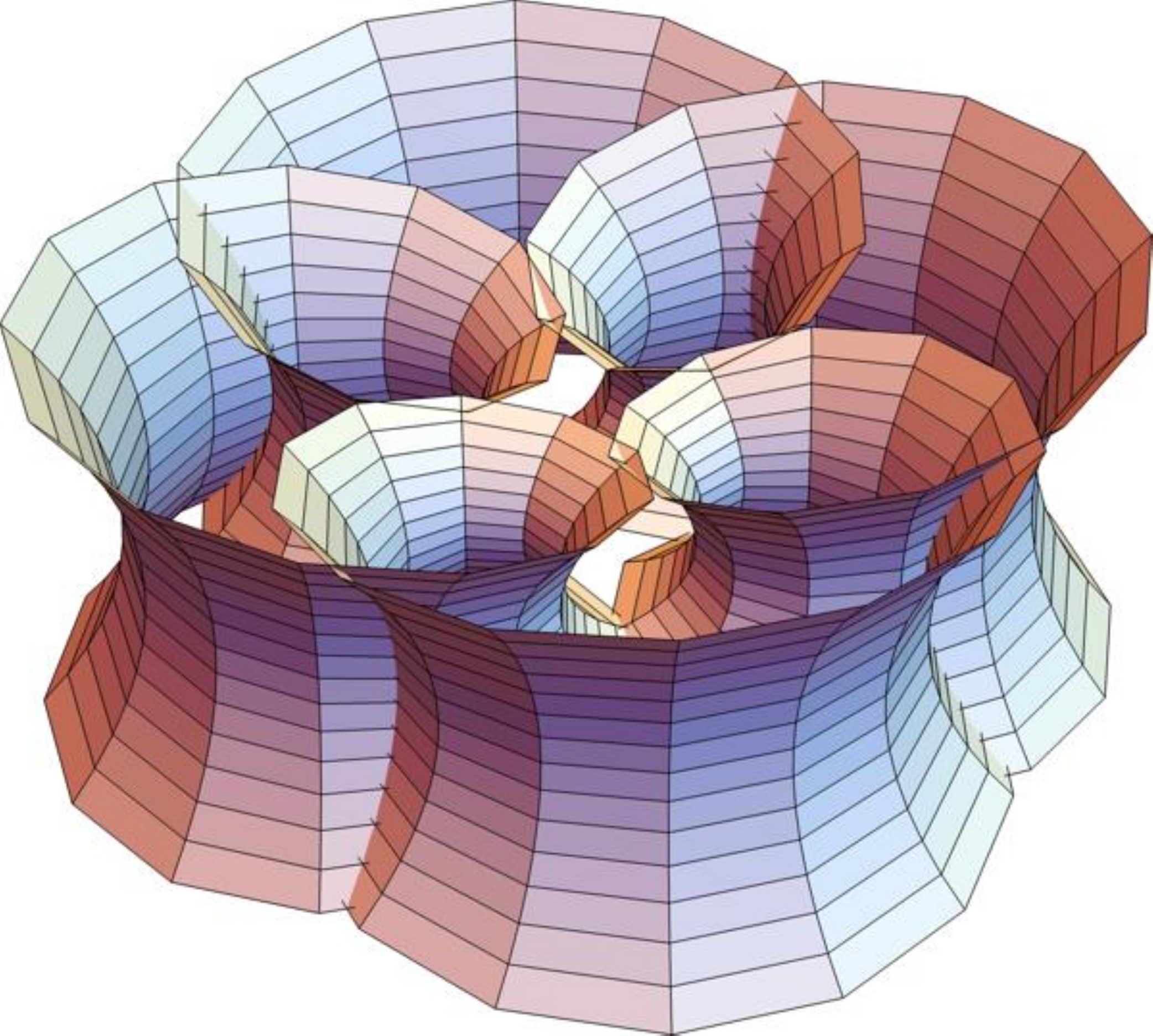}
		\end{minipage}
		
		\begin{minipage}{0.3\hsize}
		\centering
		\includegraphics[scale=0.15]{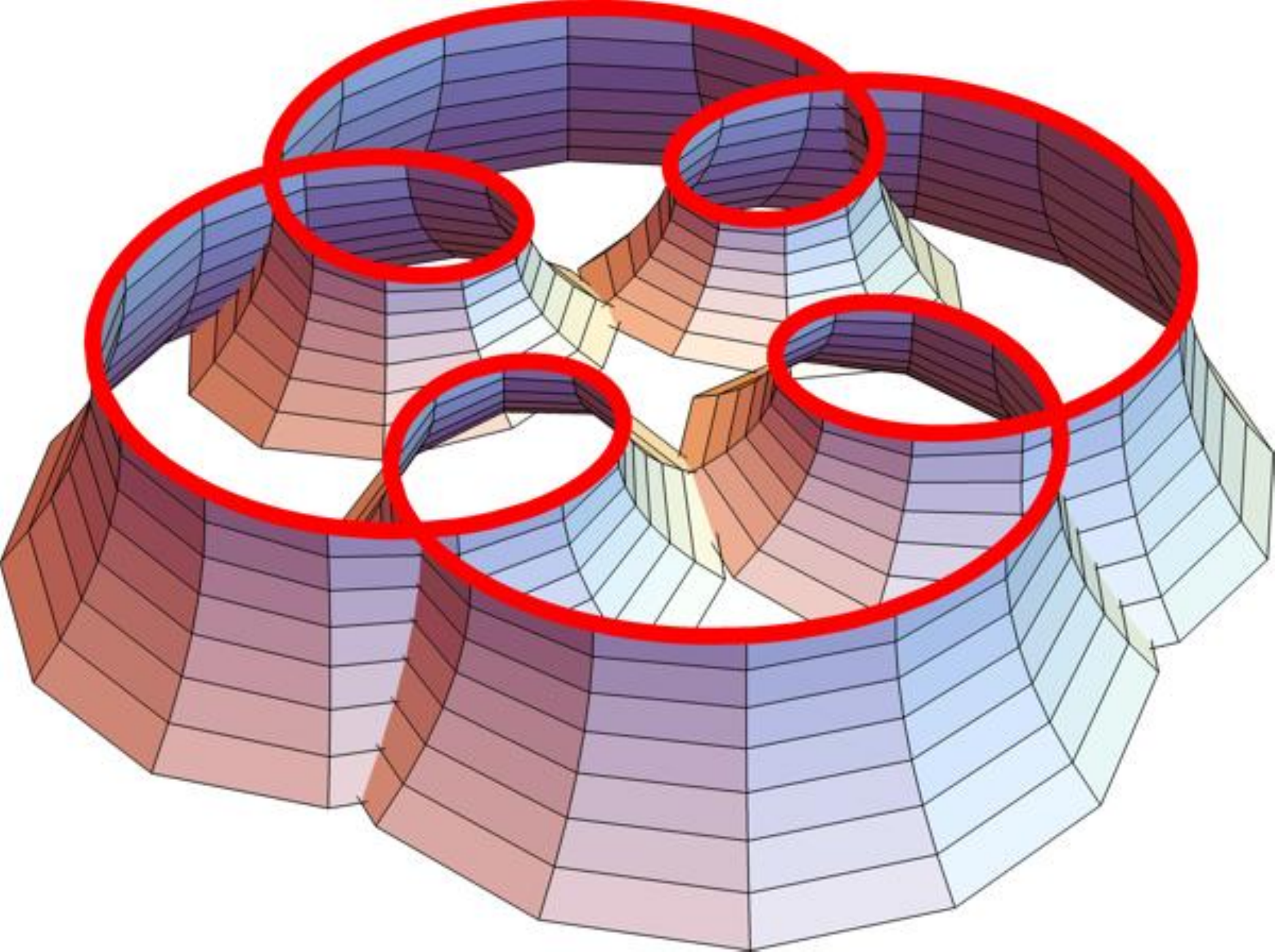}
		\end{minipage}
	\end{tabular}
\caption{The surface with $k=3$ (left), and a half cut away from the surface by the $x_1x_2$-plane (right).}
\end{figure}

\section*{Acknowledgement}
The author would like to thank {the referee, and} Professor Shoichi Fujimori {and Keisuke Teramoto} for valuable comments and suggestions.

\bigskip
{\small \noindent
Shin Kaneda \\
Department of Mathematics \\
Hiroshima University \\
Higashihiroshima, Hiroshima 739-8526, Japan \\ {\itshape E-mail address}\/: \href{mailto:shin-kaneda@hiroshima-u.ac.jp}{shin-kaneda@hiroshima-u.ac.jp}
}
\end{document}